\documentclass[titlepage,11pt]{article}
\usepackage{latexsym}
\usepackage{amsfonts}
\usepackage{amsmath}
\newcommand\blackslug{\hbox{\hskip 1pt \vrule width 4pt height 8pt depth 1.5pt
        \hskip 1pt}}
\newcommand\bbox{\hfill \quad \blackslug \bigbreak}
\def\d{\hbox{-}}

\def\l{,\ldots,}

\title{Disjoint paths in tournaments}
\author{Maria Chudnovsky\thanks{Supported by NSF grants DMS-0758364 and DMS-1001091.}\\
Columbia University, New York, NY 10027, USA
\\
\\
Alex Scott\\
Mathematical Institute,
University of Oxford,
Oxford
OX2 6GG, UK
\\
\\
Paul Seymour\thanks{Supported by NSF grant DMS-0901075 and ONR grant N00014-10-1-0680.}\\
Princeton University, Princeton, NJ 08544, USA}

\date{}

\newtheorem{thm}{}[section]

\newcommand{\Proof}{\noindent{\bf Proof.}\ \ }

\begin{document}
\maketitle
\begin{abstract}
Given $k$ pairs of vertices $(s_i,t_i)\;(1\le i\le k)$ of a digraph $G$, how can we test whether there exist $k$ vertex-disjoint directed
paths from $s_i$ to $t_i$ for $1\le i\le k$? This is NP-complete in general digraphs, even for $k = 2$~\cite{FHW},
but for $k=2$
there is a polynomial-time algorithm when $G$ is a tournament (or more generally, a semicomplete digraph), 
due to Bang-Jensen and Thomassen~\cite{bangtho}. Here we prove that for all fixed $k$ there 
is a polynomial-time algorithm to solve the problem when $G$ is semicomplete.
\end{abstract}
\section{Introduction}

Let $s_1,t_1\l s_k,t_k$ be vertices of a graph or digraph $G$. The {\em $k$ vertex-disjoint paths problem} is to determine whether
there exist vertex-disjoint paths $P_1\l P_k$ (directed paths, in the case of a digraph) 
such that $P_i$ is from $s_i$ to $t_i$ for $1\le i\le k$.
For undirected graphs, this problem is solvable in polynomial time for all fixed $k$; this was one of the highlights
of the Graph Minors project of Robertson and the third author~\cite{RS}.
The directed version is a natural and important question, but it was shown by Fortune, Hopcroft and Wyllie~\cite{FHW} that, 
without further restrictions on the input $G$, this problem is NP-complete for digraphs, even for $k = 2$. 
This motivates the study of subclasses of digraphs for which the problem is polynomial-time solvable.

In this paper, all graphs and digraphs are finite, and without loops or parallel edges; thus if $u,v$ are distinct vertices
of a digraph then there do not exist two edges both from $u$ to $v$, although there may be edges
$uv$ and $vu$. Also, by a ``path'' in a digraph we always mean a directed path.
A digraph is a {\em tournament} if for every pair of distinct vertices $u,v$, exactly one of $uv,vu$ is an edge; and a digraph is
{\em semicomplete} if for all distinct $u,v$, at least one of $uv,vu$ is an edge. It was shown by Bang-Jensen and Thomassen~\cite{bangtho}
that 
\begin{itemize}
\item the $k$ vertex-disjoint paths problem (for digraphs) is NP-complete if $k$ is not fixed, even when $G$ is a tournament; 
\item the two vertex-disjoint paths problem is solvable in polynomial time if $G$ is semicomplete.
\end{itemize}
We shall show:

\begin{thm}\label{tourthm}
For all fixed $k\ge 0$, the $k$ vertex-disjoint paths problem is solvable in polynomial time if $G$ is semicomplete.
\end{thm}

In fact we will prove a result for a wider class of digraphs, that we define next.
Let $P$ be a path of a digraph $G$, with vertices $v_1\l v_n$ in order. We say $P$ is {\em minimal} if 
$j\le i+1$ for every edge $v_iv_j$ of $G$ with $1\le i,j\le n$.
Let $d\ge 1$; we say that a digraph $G$ is {\em $d$-path-dominant} if for every minimal path $P$ of $G$ with $d$ vertices, every vertex
of $G$ either belongs to $V(P)$ or has an out-neighbour in $V(P)$ or has an in-neighbour in $V(P)$. Thus a digraph is $1$-path-dominant if and only if
it is semicomplete; and $2$-path-dominant if and only if its underlying simple graph is complete multipartite. We will show:

\begin{thm}\label{mainthm}
For all fixed $d,k\ge 1$, the $k$ vertex-disjoint paths problem is solvable in polynomial time if $G$ is $d$-path-dominant.
\end{thm}

We stress here that we are looking for vertex-disjoint paths. One can ask the same for edge-disjoint paths, and that question has also been recently solved
for tournaments, and indeed for digraphs with bounded independence number~\cite{FS}, but the solution is completely different. We do not know
a polynomial-time algorithm for the two vertex-disjoint paths problem for digraphs with independence number two. 

But we can extend \ref{mainthm} in a different way:

\begin{thm}\label{mainlength}
For all $d,k\ge 1$, there is a polynomial-time algorithm as follows:
\begin{itemize}
\item{{\bf Input:}} Vertices $s_1,t_1\l s_k,t_k$ of a $d$-path-dominant digraph $G$, and integers $x_1\l x_k\ge 1$.
\item{{\bf Output:}} Decides whether there exist pairwise
vertex-disjoint directed paths $P_1\l P_k$ of $G$ such that 
for $1\le i\le k$, $P_i$ is from $s_i$ to $t_i$ and has at most $x_i$ vertices.
\end{itemize}
\end{thm}

Let $s_1,t_1\l s_k, t_k$ be vertices of a digraph $G$. We call
$(G,s_1,t_1\l s_k,t_k)$ a {\em problem instance}.
A {\em linkage} in a digraph $G$ is a sequence $L = (P_i\;:1\le i\le k)$ of vertex-disjoint 
paths, and $L$ is a linkage {\em for} a problem instance $(G,s_1,t_1\l s_,t_k)$ if $P_i$ is from $s_i$ to $t_i$ for each $i$.
(With a slight abuse of notation, we shall call $k$ the ``cardinality'' of $L$, and $P_1\l P_k$ its ``members''. Also,
every subsequence of $(P_i\;:1\le i\le k)$ is a linkage $L'$, and we say $L$ ``includes'' $L'$. )
If $x = (x_1\l x_k)$ is a $k$-tuple of integers, we say a linkage 
$(P_i\;:1\le i\le k)$ is an {\em $x$-linkage} if each $P_i$
has $x_i$ vertices. We say a $k$-tuple of integers $x = (x_1\l x_k)$ is a {\em quality} of $(G,s_1, t_1 \l s_k, t_k)$ if 
there is an $x$-linkage for $(G,s_1,t_1\l s_,t_k)$. If $x= (x_1\l x_k)$ and $y = (y_1\l y_k)$, we say $x\le y$ if $x_i\le y_i$ for $1\le i\le k$;
and $x<y$ if $x\le y$ and $x\ne y$. We say a quality $x$ of $(G,s_1, t_1 \l s_k, t_k)$ is {\em key} 
if there is no quality $y$ with $y<x$. Our main result is the following:

\begin{thm}\label{mainchar}
For all $d,k$, there is an algorithm as follows:
\begin{itemize}
\item{{\bf Input:}} A problem instance $(G,s_1,t_1\l s_k,t_k)$ where $G$ is $d$-path-dominant.
\item{{\bf Output:}} The set of all key qualities of $(G,s_1,t_1 \l s_k,t_k)$.
\item{{\bf Running time:}} $O(n^t)$ where $t = 6k^2d(k+d) +13k$.
\end{itemize}
\end{thm}

The idea of the algorithm for \ref{mainthm} is easy described. We define an auxiliary digraph $H$ with two special vertices $s_0, t_0$, and prove that there is a 
path in $H$ from $s_0$ to $t_0$ if and only if there is a linkage for $(G,s_1,t_1\l s_,t_k)$. Thus to solve the problem of \ref{mainthm} it suffices
to construct $H$ in polynomial time. The more general question of \ref{mainchar} is solved similarly, by assigning appropriate weights to the edges of $H$.

Recently we have been able to extend \ref{tourthm} to a more general class of digraphs, namely the digraphs whose vertex set can be 
partitioned into a bounded number of subsets such that each subset induces a semicomplete digraph. The proof is by a modification
of the method of this paper, but it is considerably more difficult and not included here.

\section{A useful enumeration}

If $P$ is a path of a digraph $G$, its {\em length} is $|E(P)|$ (every path has at least one vertex); 
and $s(P), t(P)$ denote the first and last vertices of $P$, respectively.
If $F$ is a subdigraph of $G$, 
a vertex $v$ of $G\setminus V(F)$ is {\em $F$-outward} if no vertex of $F$ is adjacent
from $v$ in $G$; and {\em $F$-inward} if no vertex of $F$ is adjacent
to $v$ in $G$. If $F$ is a digraph and $v\in V(F)$, $F\setminus v$ denotes the digraph obtained from $F$ by deleting $v$; if $X\subseteq V(F)$,
$F|X$ denotes the subdigraph of $F$ induced on $X$; and $F\setminus X$ denotes the subdigraph obtained by deleting all vertices in $X$.

Now let $L = (P_i\;:1\le i\le k)$ be a linkage in $G$. We define $V(L)$ to be $V(P_1)\cup\cdots\cup V(P_k)$. A vertex $v$
is an {\em
internal vertex} of $L$ if $v\in V(L)$, and $v$ is not an end of any member of $L$. A linkage $L$ is {\em internally disjoint} from a linkage $L'$
if no internal vertex of $L$ belongs to $V(L')$ (note that this does not imply that $L'$ is internally disjoint from $L$); 
and we say that $L,L'$ are {\em internally disjoint} if each of them is internally disjoint from the other (and thus all vertices in $V(L)\cap V(L')$  must be ends of
paths in both $L$ and $L'$)

Let $Q,R$ be vertex-disjoint paths of a digraph $G$. A {\em planar $(Q,R)$-matching}
is a linkage $(M_j\;:1\le j\le n)$ for some $n\ge 0$, such that 
\begin{itemize}
\item $M_1\l M_n$ each have either two or three vertices;
\item $s(M_1)\l s(M_n)$ are vertices of $Q$, in order in $Q$; and
\item $t(M_1)\l t(M_n)$ are vertices of $R$, in order in $R$.
\end{itemize}

Fix $d,k\ge 1$, and let $L = (P_1\l P_k)$ be a linkage in a $d$-path-dominant digraph $G$. A subset $B\subseteq V(L)$ is said to be {\em acceptable}
(for $L$) if
\begin{itemize}
\item for $1\le j\le k$, if $uv$ is an edge of $P_j$ and $v\in B$ then $u\in B$
(and so $Q_j = P_j|B$ and $R_j = P_j|(V(G)\setminus B)$ are paths if they are non-null);
\item for $1\le i,j\le k$, there is no planar $(Q_i, R_j)$-matching of cardinality $(k-1)d+k^2+2$ internally disjoint from $L$.
\end{itemize}
Thus $\emptyset$ and $V(L)$ are acceptable.

\begin{thm}\label{acceptable}
Let $d\ge 1$, let $(G,s_1,t_1\l s_k,t_k)$ be a problem instance, where 
$G$ is $d$-path-dominant, let $x$ be a key quality,
and let $L = (P_1\l P_k)$ be an $x$-linkage for $(G,s_1,t_1\l s_k,t_k)$.
Suppose that $B\subseteq V(L)$ is acceptable for $L$ and $B\ne V(L)$. Then there exists $v\in V(L)\setminus B$ such that 
$B\cup \{v\}$ is acceptable for $L$.
\end{thm}
\Proof Let $A = V(G)\setminus B$. For $1\le j\le k$, let $Q_j = P_j|B$ and $R_j = P_j|A$. Let $q_j, r_j$ be the last vertex of $Q_j$ 
and the first vertex of $R_j$, respectively (if they exist).
\\
\\
(1) {\em For $1\le j\le k$, $P_j$ is a minimal path of $G$. In particular, 
the only edge of $G$ from $V(Q_j)$ to $V(R_j)$ (if there is one) is $q_jr_j$. Moreover, every three-vertex path from $V(Q_j)$
to $V(R_j)$ with internal vertex in $V(G)\setminus V(L)$ uses at least one of $q_j, r_j$. Consequently,
there is no planar $(Q_j, R_j)$-matching of cardinality three internally disjoint from $L$.}
\\
\\
For suppose there is an edge $uv$ of $G$ such that $u,v\in V(P_j)$ and $u$ is before $v$ in $P_j$, and there is at least
one vertex of $P_j$ between $u$ and $v$. If we delete from $P_j$ the vertices of $P_j$
strictly between $u$ and $v$, and add the edge $uv$, we obtain a path from $s_j$ to $t_j$ disjoint from
every member of $L$ except $P_j$, and with strictly fewer vertices than $P_j$,
contradicting that $x$ is key. Thus $P_j$ is induced.
Similarly there is no three-vertex path from $V(Q_j)$
to $V(R_j)$ with internal vertex in $V(G)\setminus V(L)$ containing neither of $q_j, r_j$. 
The final assertion follows. This proves (1).

\bigskip
From (1), the theorem holds if $k = 1$, so we may assume that $k\ge 2$.
\\
\\
(2) {\em We may assume that for all $i\in \{1\l k\}$, if $R_i$ is non-null then 
for some $j\in \{1\l k\}$ with $j\ne i$, there is a planar $(Q_i, R_j\setminus r_j)$-matching of cardinality $(k-1)d + k^2$ internally disjoint from $L$.}
\\
\\
For suppose that some $i$ does not satisfy the statement of (2). Thus $R_i$ is non-null, and 
there is no $j$ as in (2). Since $R_i$ is non-null, it follows that $r_i$ exists.
We may assume that $B\cup \{r_i\}$ is not acceptable. Consequently, one of the two conditions in the definition of
``acceptable'' is not satisfied by $B\cup \{r_i\}$. The first is satisfied since $r_i$ is the first vertex of $R_i$.
Thus the second is false, and so for some $i',j\in \{1\l k\}$, 
there is a planar $(P_{i'}|(B\cup \{r_i\}), P_{j}|(A\setminus \{r_i\}))$-matching of cardinality $(k-1)d+k^2+2$ internally disjoint from $L$.
Since there is no planar $(Q_{i'}, R_{j})$-matching of cardinality $(k-1)d+k^2+2$ internally disjoint from $L$, and $P_{j}|(A\setminus \{r_i\})$ is a subpath of $R_{j}$,
it follows that $P_{i'}|(B\cup \{r_i\})\ne Q_{i'}$, and so $i' = i$. Since only one vertex of $P_{i}|(B\cup \{r_i\})$ does not belong to $Q_i$, 
it follows that there is a planar $(Q_{i}, R_{j}\setminus r_j)$-matching of cardinality $(k-1)d+k^2$ internally disjoint from $L$.
Since $(k-1)d+k^2\ge 4$ (because $k\ge 2$), (1) implies that $j\ne i$. This proves (2).
\\
\\
(3) {\em We may assume that for some $p\ge 2$, and for all $i$ with $1\le i < p$, there is a planar $(Q_{i}, R_{i+1}\setminus r_{i+1})$-matching 
of cardinality $(k-1)d + k^2$ internally disjoint from $L$,
and there is a planar $(Q_{p}, R_{1}\setminus r_1)$-matching of cardinality $(k-1)d + k^2$ internally disjoint from $L$.}
\\
\\
For by hypothesis, there exists $i\in \{1\l k\}$ such that $R_i$ is non-null. By repeated application of (2), there exist distinct
$h_1\l h_p\in\{1\l k\}$ such that for $1\le i\le p$
there is a planar $(Q_{h_i}, R_{h_{i+1}}\setminus r_{h_{i+1}})$-matching of cardinality $(k-1)d + k^2$ internally disjoint from $L$, where $h_{p+1} = h_1$;
and $p\ge 2$ by (1). Without loss
of generality, we may assume that $h_i = i$ for $1\le i\le p$. This proves (3).

\bigskip
Let us say a planar $(Q,R)$-matching is {\em $s$-spaced} if no subpath of $Q$ with at most $s$ vertices meets more than one member of
the matching, and no subpath of $R$ with at most $s$ vertices meets more than one member of
the matching.
\\
\\
(4) {\em We may assume that for some $p\ge 2$, and for all $i$ with $1\le i < p$, there is a planar $(Q_{i}, R_{i+1}\setminus r_{i+1})$-matching $L_i$,
and there is a planar $(Q_{p}, R_{1}\setminus r_1)$-matching $L_p$, such that 
\begin{itemize}
\item $L_1\l L_p$ all have cardinality $k$;
\item they are pairwise internally disjoint;
\item each of $L_1\l L_p$ is internally disjoint from $L$; and
\item each of $L_1\l L_p$  is $(d+1)$-spaced.
\end{itemize}
}
For let $L_i'$ be a planar $(Q_{i}, R_{i+1}\setminus r_{i+1})$-matching of cardinality $(k-1)d + k^2$ internally disjoint from $L$, for $1\le i<p$, and
let $L_p'$ be a planar $(Q_{p}, R_{1}\setminus r_1)$-matching of cardinality $(k-1)d + k^2$ internally disjoint from $L$. We choose $L_i\subseteq L_{i}'$ inductively.
Suppose that for some $h<p$, we have chosen $L_1\l L_h$, such that
\begin{itemize}
\item $L_1\l L_h$ all have cardinality $k$;
\item they are pairwise internally disjoint;
\item each of $L_1\l L_h$ is internally disjoint from $L$; and
\item each of $L_1\l L_h$  is $(d+1)$-spaced.
\end{itemize}
We define $L_{h+1}$ as follows. The union of the sets of internal vertices of $L_1\l L_h$ has cardinality at most $hk\le k(k-1)$, and
so $L_{h+1}'$ includes a planar $(Q_{h+1}, R_{h+2}\setminus r_{h+2})$-matching (or $(Q_{p}, R_{1}\setminus r_1)$-matching, if $h = p-1$) 
of cardinality $(k-1)d + k^2-k(k-1)= 1+(k-1)(d+1)$, internally disjoint from each of $L_1\l L_h$.
By ordering the members of this matching in their natural order, and taking only the $i$th terms, where $i = 1, 1+ (d+1), 1+2(d+1)\ldots$, 
we obtain a $(d+1)$-spaced
matching of cardinality $k$. Let this be $L_{h+1}$. This completes the inductive definition of $L_1\l L_p$, and so proves (4).

\bigskip
For $1\le i \le p$, let $L_i = \{M^1_i\l M^k_i\}$, numbered in order; thus, if $q^h_i$ and $r^h_{i+1}$ denote the first and last vertices of $M^h_i$,
then $q^1_i\l q^k_i$ are distinct and in order in $Q_i$, and $r_{i+1}, r^1_{i+1}\l r^k_{i+1}$ are
distinct and in order in $R_{i+1}$ (or in $R_1$ if $i = p$). For $1\le i\le p$ and $2\le h\le k$, let $Q^h_i$ be the subpath of $P_i$
with $d$ vertices and with last vertex $q^h_i$. (Thus $q^{h-1}_i$ does not belong to $Q^h_i$ since $L_i$ is $d$-spaced, and indeed $(d+1)$-spaced.)
Since $P_i$ and hence $Q^h_i$ is a minimal path of $G$, and $G$ is $d$-path-dominant, it follows that for $1\le i\le p$ and $2\le h\le k$,
$r^{h-1}_i$ is adjacent to or from some vertex $v$ of $Q^h_i$. Since $r^{h-1}_i\ne r_i$, (1) implies that $r^{h-1}_i$ is not adjacent from any
vertex of $Q^h_i$; and so there is a path $R^{h-1}_i$ from $r^{h-1}_i$ to $q^h_i$ of length at most $d$, such that all its internal
vertices belong to $Q^h_i$. For $1\le i\le p$, and $1\le h<k$, let $S^h_i$ be the path
$$q^h_i\d M^h_i\d r^h_{i+1}\d R^h_{i+1} \d q^{h+1}_{i+1},$$
or 
$$q^h_p\d M^h_p\d r^h_{1}\d R^h_{1}\d q^{h+1}_1$$
if $i = p$; then $S^h_i$ is a path from $q^h_i$ to $q^{h+1}_{i+1}$ (or to $q^{h+1}_{1}$ if $i = p$), of length at most $d+2$.
Thus (reading subscripts modulo $p$) concatenating 
$S^1_i, S^2_{i+1}\l S^{p-1}_{i+p-2}$ and $M^p_{i-1}$ gives a path $T_i'$ from $q^1_i$ to $r^p_{i}$ of length at most
$(p-1)(d+2) + 2$. The subpath $T_i$ of $P_i$ from $q^1_i$ to $r^p_{i}$ has length at least $(p+k-2)(d+1)+2$, since $L_{i-1}, L_i$ are $(d+1)$-spaced and 
$r_i$ is different from $r^1_i$; and since $p+k-2\ge 2(p-1)$ and $d+1> (d+2)/2$, it follows that $T_i$ has length strictly 
greater than that of $T_i'$.  Let $P_i'$ be obtained from $P_i$ by replacing the
subpath $T_i$ by $T_i'$, for $1\le i\le p$, and let $P_{i'} = P_i$ for $p+1\le i\le k$.
Then $\{P_1'\l P_k'\}$ is a linkage for $(G,s_1,t_1\l s_,t_k)$, contradicting that $x$ is key.
This proves \ref{acceptable}.~\bbox

We deduce:

\begin{thm}\label{vertexorder}
Let $d\ge 1$, let $(G,s_1,t_1\l s_k,t_k)$ be a problem instance where $G$ is $d$-path-dominant, let $x$ be a key quality,
and let $L = (P_1\l P_k)$ be an $x$-linkage for $(G,s_1,t_1\l s_k,t_k)$. Let $c = (k-1)d+k^2+2$.
Then there is an enumeration $(v_1\l v_n)$ of $V(L)$, such that 
\begin{itemize}
\item for $1\le h\le k$ and $1\le p,q\le n$, if $v_pv_q$ is an edge of $P_h$ then $p<q$;
\item for $1\le h,i\le k$ and $1\le p\le n-1$, and every $cd$-vertex subpath $Q$ of $P_h|\{v_1\l v_p\}$, and every $cd$-vertex subpath
$R$ of $P_i|\{v_{p+1}\l v_n\}$, there are at most $c(2k+1)$ vertices of $G$ that are both $Q$-outward and $R$-inward.
\end{itemize}
\end{thm}
\Proof
Since $\emptyset$ is acceptable for $L$, by repeated application of \ref{acceptable} 
implies that there is an enumeration $(v_1\l v_n)$ of $V(L)$, such that
$\{v_1\l v_p\}$ is acceptable for $0\le p\le n$. We claim that this enumeration satisfies the theorem. For certainly the first bullet
holds; we must check the second. Thus, let $1\le p\le n$, and let $B = \{v_1\l v_p\}$ and $A = \{v_{p+1}\l v_n\}$. For
$1\le h\le k$, let $Q_h = P_h|B$ and $R_h = P_h|A$. Now let $1\le h,i\le k$, and let $Q, R$ be $cd$-vertex subpaths of $Q_h, R_i$ respectively.
Let  $X$ be the set of all vertices of $G$ that are both $Q$-outward and $R$-inward. We must show that $|X|\le c(2k+1)$.
\\
\\
(1) {\em If $x_1\l x_c\in X$ are distinct, then there exist $y_1\l y_c\in V(Q)$, distinct and in order in $Q$, such that
$y_jx_j$ is an edge for $1\le j\le c$.}
\\
\\
For $Q$ has $cd$ vertices; let its vertices be $q_1\l q_{cd}$ in order. Let $1\le j\le c$. The subpath of $Q$ induced on $\{q_s\::(j-1)d<s\le jd\}$
has $d$ vertices, and since $Q$ is a minimal path of $G$ and $G$ is $d$-path-dominant, and $X\cap V(Q) = \emptyset$, it follows that
$x_j$ is in- or out-adjacent to a vertex of this subpath, say $y_j$. Since
$x_j\in X$ and hence is $Q$-outwards, it follows that $x_jy_j$ is not an edge, and so $y_jx_j$ is an edge. But then $y_1\l y_c$
satisfy (1). This proves (1).
\\
\\
(2) {\em The sets $X\setminus V(L)$, $X\cap V(Q_g)\; (1\le g\le k)$ and $X\cap V(R_g)\;(1\le g\le k)$ all have cardinality at most $c-1$, and
hence $|X|\le (2k+1)(c-1)$.}
\\
\\
For suppose that there exist distinct $x_1\l x_c\in X\setminus V(L)$. By (1) there exist distinct $y_1\l y_c\in V(Q)$, in order in $Q$, such that
$y_jx_j$ is an edge for $1\le j\le c$; and similarly there exist $z_1\l z_c\in V(R)$, in order in $R$, such that
$x_jz_j$ is an edge for $1\le j\le c$. But then the $c$ paths $y_j\d x_j\d z_j\;(1\le j\le c)$ form a planar $(Q_h,R_i)$-matching
of cardinality $c$, internally disjoint from $L$, contradicting that $\{v_1\l v_p\}$ is acceptable. Thus
$|X\setminus V(L)|\le c-1$. Now suppose that for some $g\in \{1\l k\}$, there exist distinct $x_1\l x_c$ in $X\cap V(R_g)$,
numbered in order in $R_g$. Choose $y_1\l y_c$ as in (1); then the paths $y_jx_j\;(1\le j\le c)$ form 
a planar $(Q_h,R_g)$-matching 
of cardinality $c$, internally disjoint from $L$, contradicting that $\{v_1\l v_p\}$ is acceptable. Thus
$|X\cap V(R_g)|\le c-1$, and similarly $|X\cap V(Q_g)|\le c-1$, for $1\le g\le k$. This proves (2).

\bigskip
From (2), the theorem follows.~\bbox

\section{Confusion and the auxiliary digraph}

Let $(G, s_1, t_1\l s_k, t_k)$ be a problem instance, and let $L = (M_1\l M_k)$ be a linkage in $G$ (not necessarily a linkage for $(G, s_1, t_1\l s_k, t_k)$).
Let 
$A(L)$ be the set of all vertices in  $V(G)\setminus V(L)$ that are $M_j\setminus t(M_j)$-inward for some $j\in \{1\l k\}$ such that $t(M_j)\ne t_j$
and let
$B(L)$ be the set of all vertices in $V(G)\setminus V(L)$ that are $M_j\setminus s(M_j)$-outward for some $j\in \{1\l k\}$ such that $s(M_j)\ne s_j$.
We call $|A(L)\cap B(L)|$ the {\em confusion} of $L$ ; and it is helpful to keep the confusion small, as we shall see.

A {\em $(k,m,c)$-rail} in a problem instance $(G, s_1, t_1\l s_k, t_k)$ is a triple $(L,X,Y)$, where
\begin{itemize}
\item $L$ is a linkage in $G$ consisting of $k$ paths $(M_1\l M_k)$ (but not necessarily a linkage for $(G, s_1, t_1\l s_k, t_k)$);
\item for $1\le j\le k$, $M_j$ has at most $2m$ vertices, and if it has fewer than $2m$ vertices then $M_j$ either has first vertex $s_j$ or last vertex $t_j$;
\item $L$ has confusion at most $c$;
\item $X,Y$ are disjoint subsets of $V(G)\setminus V(L)$; and
\item $X\subseteq A(L)$, $Y\subseteq B(L)$, and $X\cup Y = A(L)\cup B(L)$.
\end{itemize}

\begin{thm}\label{countrails}
For all $k,m,c\ge 0$, if $(G, s_1, t_1\l s_k, t_k)$ is a problem instance and $G$ has $n$ vertices then
there are at most $2^cn^{2km}(2km)^k$ $(k,m,c)$-rails in $(G, s_1, t_1\l s_k, t_k)$.
Moreover, for all fixed $k,m,c\ge 0$, there is an algorithm which, with input a problem instance $(G, s_1, t_1\l s_k, t_k)$, 
finds all its $(k,m,c)$-rails in time $O(n^{2km+1})$, where $n = |V(G)|$.
\end{thm}
\Proof
First, if $L$ is a linkage with $k$ paths each with at most $2m$ vertices, then $|V(L)|\le 2km$,
and so the number of such linkages is at most $n^{2km}(2km)^k$, as is easily seen. Now fix a linkage $L$ satisfying the first two bullets
in the definition of $(k,m,c)$-rail; let us count the number
of pairs $(X,Y)$ such that $(L,X,Y)$ is a $(k,m,c)$-rail. There are none unless $|A(L)\cap B(L)|\le c$;
and in that case, there are at most $2^c$ possibilities for the pair $(X,Y)$,
since $X$ consists of $A(L)\setminus B(L)$ together with some subset of $A(L)\cap B(L)$, and
$Y = (A(L)\cup B(L)) \setminus X$.

For the algorithm, we first find all linkages $L$ with $k$ paths each with at most $2m$ vertices, by examining all ordered $2km$-tuples of distinct vertices
of $G$. For each such $L$, we check whether it satisfies the first three bullets in the definition of $(k,m,c)$-rail (this takes time $O(n)$);
if not we discard it and otherwise we partition
$A(L)\cap B(L)$ into two subsets in all possible ways, and output the corresponding $(k,m,c)$-rails.  
The result follows.~\bbox

Let $(L,X,Y)$ and $(L', X', Y')$ be distinct $(k,m,c)$-rails in $G$, and let $L = (P_1\l P_k)$ and $L' = (P_1'\l P_k')$.
We write $(L,X,Y)\rightarrow (L', X', Y')$ if the following hold:
\begin{itemize}
\item for $1\le i\le k$, $P_i\cup P_i'$ is a path from the first vertex of $P_i$ to the last vertex of $P_i'$;
\item for $1\le i\le k$, $V(P_i')\subseteq V(P_i)\cup X$, and $V(P_i)\subseteq V(P_i')\cup Y'$; and
\item $X'\subseteq X$, and $Y\subseteq Y'$.
\end{itemize}

Let $(G, s_1, t_1\l s_k, t_k)$ be a problem instance, and let $\mathcal{T}$  be the set of all $(k,m,c)$-rails in 
$(G, s_1, t_1\l s_k, t_k)$. Take two new vertices $s_0, t_0$, and let us define a digraph $H$ with vertex set $\mathcal{T}\cup \{s_0, t_0\}$
as follows. Let $u,v\in V(H)$. If $u,v\in \mathcal{T}$ are distinct, then $uv\in E(H)$ if and only if $u\rightarrow v$.
If $u = s_0$ and $v\in \mathcal{T}$, let $v = (L,X,Y)$ where $L = (M_1\l M_k)$; then $uv\in E(H)$ if and only if $M_j$ has first vertex $s_j$ for 
all $j\in \{1\l k\}$.
Similarly, if $u \in  \mathcal{T}$ and $v= t_0$, let $u = (L,X,Y)$ where $L = (M_1\l M_k)$; 
then $uv\in E(H)$ if and only if $M_j$ has last vertex $t_j$ for all $j\in \{1\l k\}$. This defines $H$. We call 
$H$ the {\em $(k,m,c)$-tracker} of $(G, s_1, t_1\l s_k, t_k)$.

We shall show that with an appropriate choice of $m,c$, when $G$ is $d$-path-dominant
we can reduce our problems about linkages for $(G, s_1, t_1\l s_k, t_k)$ to problems about
paths from $s_0$ to $t_0$ in the $(k,m,c)$-tracker.
Let $(G, s_1, t_1\l s_k, t_k)$ be a problem instance, let $(P_1\l P_k)$ be a linkage for  $(G, s_1, t_1\l s_k, t_k)$, and let
$P$ be a path from $s_0$ to $t_0$ in the $(k,m,c)$-tracker. Let $P$ have vertices 
$$s_0, (L_1,X_1,Y_1)\l (L_n,X_n,Y_n), t_0$$
in order, and 
let $L_p = (M_{p,1}\l M_{p,k})$ for $1\le p\le n$. We say that $P$ {\em traces} $(P_1\l P_k)$ if 
$P_j$ is the union of $M_{1,j}\l M_{n,j}$ for all $j\in \{1\l k\}$.

\begin{thm}\label{pathtolink}
Let $k,m,c\ge 0$ be integers, and let $(G, s_1, t_1\l s_k, t_k)$ be a problem instance, with $(k,m,c)$-tracker $H$. Every path in $H$ from $s_0$ to $t_0$
traces some linkage for $(G, s_1, t_1\l s_k, t_k)$.
\end{thm}
\Proof Let $P$ be a path of $H$, with vertices 
$$s_0, (L_1,X_1,Y_1)\l (L_n,X_n,Y_n), t_0$$
in order, and  
let $L_p = (M_{p,1}\l M_{p,k})$ for $1\le p\le n$. For $1\le p\le n$ and $1\le j\le k$, let $P_{p,j}$ be the union of
$M_{1,j}\l M_{p,j}$. 
\\
\\
(1) {\em For $1\le p\le n$ and $1\le j \le k$, every vertex of $P_{p,j}$ belongs to $Y_p\cup V(M_{p,j})$.}
\\
\\
We prove this by induction on $p$. If $p = 1$ the claim is true, since then $P_{1,j} = M_{1,j}$. We assume then that $p>1$ and the result holds
for $p-1$. Let $v\in V(P_{p,j})$. If $v\in V(M_{p,j})$ then the claim is true, so we assume not. Since $v\in V(P_{p,j})$, and $P_{p,j} = P_{p-1,j}\cup M_{p,j}$,
it follows that $v\in V(P_{p-1,j})$, and so from the inductive hypothesis, $v\in Y_{p-1}\cup V(M_{p-1,j})$. But since
$(L_{p-1},X_{p-1},Y_{p-1})\rightarrow (L_p, X_p, Y_p)$, we deduce that $Y_{p-1}\subseteq Y_p$, and $V(M_{p-1,j})\subseteq V(M_{p,j})\cup Y_p$, and so 
$v\in V(M_{p,j})\cup Y_p$. This proves (1).
\\
\\
(2) {\em For $1\le p\le n$ and $1\le j \le k$, $P_{p,j}$ is a path from $s_j$ to the last vertex of $M_{p,j}$.}
\\
\\
The claim holds if $p = 1$; so we assume that $p>1$ and the claim holds for $p-1$. Thus $P_{p-1,j}$ is a path from $s_j$ to the last vertex of $M_{p-1,j}$; 
and also, $M_{p-1,j}\cup M_{p,j}$ is a path,
from the first vertex of $M_{p-1,j}$ to the last vertex of $M_{p,j}$, since $(L_{p-1},X_{p-1},Y_{p-1})\rightarrow (L_p, X_p, Y_p)$.
We claim that every vertex $v$ that belongs to both of $P_{p-1,j}, M_{p,j}$ also belongs to $M_{p-1,j}$. For suppose not; then by (1),
$v\in Y_{p-1}$ since $v\in V(P_{p-1,j})\setminus V(M_{p-1,j})$, and $v\in X_{p-1}$, since $v\in V(M_{p,j})\setminus V(M_{p-1,j})$. This is impossible since
$X_{p-1}\cap Y_{p-1} = \emptyset$. This proves that every vertex that belongs to both of $P_{p-1,j}, M_{p,j}$ also belongs to $M_{p-1,j}$. Since
$M_{p-1,j}$ is non-null, we deduce that $P_{p-1,j} \cup M_{p,j}$ is a path from $s_j$ to the last vertex of $M_{p,j}$. This proves (2).
\\
\\
(3) {\em For $1\le p\le n$, the paths $P_{p,1}\l P_{p,k}$ are pairwise vertex-disjoint.}
\\
\\
For again we proceed by induction on $p$, and may assume that $p>1$ and the result holds for $p-1$. Suppose that $v$ belongs to two of the paths
$P_{p,1}\l P_{p,k}$, say to $P_{p,1}$ and $P_{p,2}$. From the inductive hypothesis, $v$ does not belong to both of $P_{p-1,1}$ and $P_{p-1,2}$, so
we may assume that $v\in V(M_{p,1})$. Now $v\notin V(M_{p,2})$, because $L_p$ is a linkage, and so $v\in V(P_{p-1,2})$. From
(1) we deduce that $v\in Y_{p-1}\cup V(M_{p-1,2})$. But $Y_{p-1}\subseteq Y_p$, and $V(M_{p-1,2})\setminus V(M_{p,2})\subseteq Y_p$, and so
$v\in Y_p$; but $Y_p\cap V(L_p) = \emptyset$ since $(L_p,X_p,Y_p)$ is a $(k,m,c)$-rail, a contradiction. This proves (3).

\bigskip
From (2) and (3) we deduce that $(P_{n,1}\l P_{n,k})$ is a linkage $L$ for $(G,s_1,t_1\l s_k, t_k)$. Thus $P$ traces $L$. This proves \ref{pathtolink}.~\bbox

The next result is a kind of partial converse; but we have to choose $m,c$ carefully, and we need $G$ to be $d$-path-dominant, 
and the proof only works for linkages that realize a key quality.

\begin{thm}\label{linktopath}
Let $d,k\ge 1$ be integers, and let
\begin{eqnarray*}
c &=& ((k-1)d+k^2+2)(2k+1)k^2\\
m &=& ((k-1)d+k^2+2)d+1.
\end{eqnarray*}
Let $(G,s_1,t_1\l s_k,t_k)$ be a problem instance where $G$ is $d$-path-dominant, let $x$ be a key quality,
and let $(P_1\l P_k)$ be an $x$-linkage for $(G,s_1,t_1\l s_k,t_k)$. 
Let $H$ be the $(k,m,c)$-tracker of $(G,s_1,t_1\l s_k,t_k)$. Then there is a path in $H$ from $s_0$ to $t_0$ tracing $(P_1\l P_k)$.
\end{thm}
\Proof
Let $L = (P_1\l P_k)$.
By \ref{vertexorder}, there is an enumeration $(v_1\l v_n)$ of $V(L)$, such that
\begin{itemize}
\item for $1\le j\le k$ and $1\le p,q\le n$, if $v_pv_q$ is an edge of $P_j$ then $p<q$;
\item for $1\le i,j\le k$ and $1\le p\le n-1$, and every $(m-1)$-vertex subpath $Q$ of $P_i|\{v_1\l v_p\}$, and every $(m-1)$-vertex subpath
$R$ of $P_j|\{v_{p+1}\l v_n\}$, there are at most $((k-1)d+k^2+2)(2k+1)$ vertices of $G$ that are both $Q$-outward and $R$-inward.
\end{itemize}
For each $v\in V(L)$, let $\phi(v) = i$ where $v = v_i$; thus $\phi$ is a bijection from $V(L)$ onto $\{1\l n\}$.

For all $p\in \{0\l n\}$ and all $j\in \{1\l k\}$, if  $\phi(s_j)\le p$, let $Q_{p,j}$ be the maximal subpath of $P_j$
with at most $m$ vertices and with last vertex $v_q$, where $q\le p$ is maximum such that $v_q\in V(P_j)$.
If $\phi(s_j)>p$, let $Q_{p,j}$ be the null digraph. Similarly, if $\phi(t_j)>p$, let
$R_{p,j}$ be the  maximal subpath of $P_j$
with at most $m$ vertices and with first vertex $v_r$, where $r>p$ is minimum such that $v_r\in V(P_j)$.
If $\phi(t_j)\le p$, let $R_{p,j}$ be the null digraph. Thus, if $Q_{p,j}, R_{p,j}$ are both non-null, then $t(Q_{p,j})$ and $s(R_{p,j})$ are
consecutive in $P_j$.

For all $p\in \{0\l n\}$ and all $j\in \{1\l k\}$, let $M_{p,j}$ be the subpath of $P_j$ defined as follows: if both $Q_{p,j},R_{p,j}$ are non-null,
$M_{p,j}$  consists of $Q_{p,j}\cup R_{p,j}$ together with the edge of $P_j$
from $t(Q_{p,j})$ to $s(R_{p,j})$, while if one of $Q_{p,j},R_{p,j}$ is null, $M_{p,j}$ equals the other (not both can be null).
We see that, for all $p,j$, $M_{p,j}$ has at most $2m$ vertices; and either it has exactly $2m$, or its first vertex is $s_j$, or its last vertex is $t_j$.
For all $p\in \{0\l n\}$, let $L_p$ be the linkage $(M_{p,1}\l M_{p,k})$.
\\
\\
(1) {\em For all $p\in \{0\l n\}$, $L_p$ has confusion at most $c$.}
\\
\\
Let $v\in A(L_p)\cap B(L_p)$, where $A(L_p), B(L_p)$ are as in the definition of confusion. Thus there exists $j\in \{1\l k\}$ such that
$v$ is $M_{p,j}\setminus t(M_{p,j})$-inward and $t(M_{p,j})\ne t_j$. Since $t(M_{p,j})\ne t_j$, it follows from the choice of $R_{p,j}$ that 
$R_{p,j}$ has exactly $m$ vertices. Moreover, $v$ is $R_{p,j}\setminus t(R_{p,j})$-inward, since $v$ is $M_{p,j}\setminus t(M_{p,j})$-inward. Similarly,
there exists $i\in \{1\l k\}$ such that $v$ is $Q_{p,i}\setminus s(Q_{p,i})$-outward and 
$Q_{p,i}$ has $m$ vertices. For each choice of $i,j\in \{1\l k\}$, there are at
most $((k-1)d+k^2+2)(2k+1)$ vertices that are both $Q_{p,i}\setminus s(Q_{p,i})$-outward and $R_{p,j}\setminus t(R_{p,j})$-inward, from the choice of the enumeration
$(v_1\l v_n)$. Consequently in total there are only $c$ possibilities for $v$, and so $|A(L_p)\cap B(L_p)|\le c$. This proves (1).
\\
\\
(2) {\em For $0\le p\le n$ and each $v\in V(L)\setminus V(L_p)$, if $\phi(v)>p$ then $v\in A(L_p)$, and if $\phi(v)\le p$ then $v\in B(L_p)$.}
\\
\\
For let $v\in V(P_j)$ say. Assume first that $\phi(v)>p$. Since $v\notin V(L_p)$, it follows that $M_{p,j}$ does not have last
vertex $t_j$; and 
since $x$ is key, $v$ is not adjacent from any vertex in $M_{p,j}$ except possibly $t(M_{p,j})$. Consequently 
$v$ is $M_{p,j}\setminus t(M_{p,j})$-inward, and hence belongs to $A(L_p)$. Similarly, if $\phi(v)\le p$ then $v\in B(L_p)$.
This proves (2).

\bigskip

For all $p\in \{0\l n\}$, define $X_p, Y_p$ as follows: 
\begin{eqnarray*}
X_p &=& \{v\in V(L)\setminus V(L_p)\::\phi(v)> p\}\cup (A(L_p)\setminus B(L_p))\\
Y_p &=& (A(L_p)\cup B(L_p))\setminus X_p.
\end{eqnarray*}
(3) {\em For all $p\in \{0\l n\}$, $(L_p, X_p, Y_p)$ is a $(k,m,c)$-rail.}
\\
\\
From (1), it suffices to check that 
\begin{itemize}
\item $X_p,Y_p$ are disjoint subsets of $V(G)\setminus V(L_p)$;
\item $X_p\subseteq A(L_p)$, $Y_p\subseteq B(L_p)$; and 
\item $X_p\cup Y_p = A(L_p)\cup B(L_p)$.
\end{itemize}
Certainly they are disjoint, and have union $A(L_p)\cup B(L_p)$. Moreover, from (2), $X_p\subseteq A(L_p)$. It remains
to show that $Y_p\subseteq B(L_p)$. Let $v\in Y_p$. Thus $v\in A(L_p)\cup B(L_p)$; and $v\notin A(L_p)\setminus B(L_p)$, since $v\notin X_p$.
Consequently $v\in B(L_p)$ as required. This proves (3).
\\
\\
(4) {\em For all $p\in \{0\l n-1\}$, and all $j\in \{1\l k\}$, 
$M_{p,j}\cup M_{p+1,j}$ is a path from the first vertex of $M_{p,j}$ to the last vertex of $M_{p+1,j}$.}
\\
\\
For $M_{p,j}, M_{p+1,j}$ are both subpaths of $P_j$, and we may assume they are distinct, and so
$v_{p+1}\in V(P_j)$. Hence, since $m>0$,  $v_{p+1}$ is the first vertex of $R_{p,j}$, and the last vertex of $Q_{p+1,j}$; and so $M_{p,j}\cup M_{p+1,j}$
is a path. Moreover, it follows from the definition of the paths $M_{p,j}$ that
$M_{p,j}\cup M_{p+1,j}$ is a path from the first vertex of $M_{p,j}$ to the last vertex of $M_{p+1,j}$. This proves (4).
\\
\\
(5) {\em For all $p\in \{0\l n-1\}$, and all $j\in \{1\l k\}$, $A(L_{p+1})\subseteq A(L_p)\cup V(L)$ and $B(L_p)\subseteq B(L_{p+1})\cup V(L)$.}
\\
\\
For let $v\in A(L_{p+1})$. We need to prove that $v\in A(L_p)\cup V(L)$, and so we may assume that $v\notin V(L)$. Choose
$j$ with $1\le j\le k$ such that
$v$ is $M_{p+1,j}\setminus t(M_{p+1,j})$-inward and $t(M_{p+1,j})\ne t_j$. Consequently $t(M_{p,j})\ne t_j$, and so if 
$v$ is $M_{p,j}\setminus t(M_{p,j})$-inward then $v\in A(L_p)$ as required, so we may assume that $v$ is adjacent from some vertex of $M_{p,j}$. 
In particular, $M_{p,j}\ne M_{p+1,j}$ and so $v_{p+1}\in V(P_j)$, and $v_{p+1} = s(R_{p,j}) = t(Q_{p+1,j})$.
Moreover, since $s(M_{p,j})$ is the only vertex of $M_{p,j}$ that may not belong to $ M_{p+1,j}$, we deduce that
$s(M_{p,j})$ is adjacent to $v$, and $s(M_{p,j})$ does not belong to $ M_{p+1,j}$.
Consequently 
$s(M_{p+1,j})\ne s_j$, and so $Q_{p+1,j}$ has $m$ vertices. Since $v$ is $M_{p+1,j}\setminus t(M_{p+1,j})$-inward, and $G$ is $d$-path-dominant,
and $M_{p+1,j}\setminus t(M_{p+1,j})$ is a minimal path of $G$, and it has $m-1\ge d+2$ vertices, 
there is a subpath of $M_{p+1,j}\setminus t(M_{p+1,j})$ with $d$ vertices, not containing
the first or second vertex of $M_{p+1,j}\setminus t(M_{p+1,j})$; and so $v$ is adjacent to some vertex $w$ of $M_{p+1,j}\setminus t(M_{p+1,j})$ different from
its first and second vertices. But $v$ is adjacent from $u$, so by replacing the subpath of $P_j$ between $u$ and $w$ by the path $u\d v\d w$, we contradict
that $x$ is key. This proves that $v\in A(L_p)$, and so $A(L_{p+1})\subseteq A(L_p)\cup V(L)$. Similarly
$B(L_p)\subseteq B(L_{p+1})\cup V(L)$. This proves (5).
\\
\\
(6) {\em For all $p\in \{0\l n-1\}$, $X_{p+1}\subseteq X_p$, and $Y_p\subseteq Y_{p+1}$.}
\\
\\
Let $v\in X_{p+1}$. Suppose first that $v\notin V(L)$. Then $v\in A(L_{p+1})\setminus B(L_{p+1})$. By (5), $v\in A(L_p)\setminus B(L_p)$, and so
$v\in X_p$ as required. Thus we may assume that $v\in V(L)$. Since $v\in X_{p+1}$, it follows that either $\phi(v)>p+1$, or $v\notin B(L_{p+1})$.
If $\phi(v)>p+1$, then
since $v\notin V(L_{p+1})$, it follows that $v\notin V(L_p)$, and hence $v\in X_p$ from the definition of $X_p$. Thus we may assume
that $\phi(v)\le p+1$ and $v\notin B(L_{p+1})$, contrary to (2). This proves that $X_{p+1}\subseteq X_p$.

For the second inclusion, let $v\in Y_p$. Suppose first that $v\notin V(L)$. Then $v\in B(L_p)$; and so $v\in B(L_{p+1})$ by (5),
and hence $v\in Y_{p+1}$ as required. Thus we may assume that $v\in V(L)$. Since $v\in Y_p$, it follows that $\phi(v) \le p$.
Now $v\notin V(L_p)$, and therefore $v\notin V(L_{p+1})$.
But $\phi(v)\le p+1$, and so by (2), $v\in B(L_{p+1})$, and consequently $v\notin X_{p+1}$. Thus $v\in Y_{p+1}$, as required. This proves that
$Y_p\subseteq Y_{p+1}$, and so proves (6).
\\
\\
(7) {\em For all $p\in \{0\l n-1\}$, and all $j\in \{1\l k\}$, 
$V(P_{p+1,j})\subseteq V(P_{p,j})\cup X_p$, and $V(P_{p,j})\subseteq V(P_{p+1,j})\cup Y_{p+1}$.}
\\
\\
To prove the first assertion,  let $v\in V(P_{p+1,j})\setminus V(P_{p,j})$. It follows that $\phi(v)>p$; but then $v\in X_p$
from the definition of $X_p$. For the second assertion, let $v\in V(P_{p,j})\setminus V(P_{p+1,j})$; then 
$\phi(v)\le p+1$, and so $v\in B(L_{p+1})$ by (2). Consequently $v\notin X_{p+1}$, and so $v\in Y_{p+1}$ as required. This proves (7).
\\
\\
(8) {\em For all $p\in \{0\l n-1\}$, $(L_p, X_p, Y_p)\rightarrow (L_{p+1}, X_{p+1}, Y_{p+1})$.}
\\
\\
This is immediate from (4), (6) and (7).

\bigskip

Now $(L_1, X_1, Y_1)\l (L_n, X_n, Y_n)$ are not necessarily all distinct. But we have:
\\
\\
(9) {\em For all $p,r$ with $0\le p\le r\le n$, if $(L_p, X_p, Y_p) = (L_r, X_r, Y_r)$, 
then $(L_p, X_p, Y_p) = (L_q, X_q, Y_q)$ for all $q$ with $p\le q\le r$.}
\\
\\
For (6) implies that $X_q\subseteq X_p$, and $X_r\subseteq X_q$, and so $X_p = X_q$, and similarly $Y_p = Y_q$. If some vertex $v$ belongs to $V(L_q)\setminus V(L_p)$,
then by (7) and (6), $v\in X_p = X_q$, a contradiction. Similarly, if $v\in V(L_p)\setminus V(L_q)$ then $v\in Y_q = Y_p$, a contradiction. This proves (9).
\\
\\
(10) {\em For all $j\in \{1\l k\}$, $M_{0,j}$ has first vertex $s_j$, and $M_{n,j}$ has last vertex $t_j$.}
\\
\\
This follows from the definitions of $M_{0,j}$ and $M_{n,j}$. 

\bigskip

We recall that $H$ is the $(k,m,c)$-tracker, with two special vertices $s_0, t_0$. 
Now (10) implies that $s_0$ is adjacent to $(L_1, X_1, Y_1)$ in $H$, and $(L_n, X_n, Y_n)$
is adjacent to $t_0$. From (8) and (9), there is a subsequence of the sequence 
$$s_0, (L_1, X_1, Y_1)\l (L_n, X_n, Y_n), t_0,$$
which lists the vertex set in order of a path of $H$ from $s_0$ to $t_0$. By \ref{pathtolink}, this path traces some linkage $L'$ for $(G, s_1, t_1\l s_k, t_k)$.
But for all $j\in \{1\l k\}$, $M_{0,j}, M_{1,j}\l M_{n,j}$ are all subpaths of $P_j$; and since their union is a path from $s_j$ to $t_j$, it follows
that their union is $P_j$. Hence $L' = L$. This proves \ref{linktopath}.~\bbox

\section{The algorithm}

Next, we need a polynomial algorithm to solve a kind of vector-valued shortest path problem.
If $n\ge 0$ is an integer, $K_n$ denotes the set of
all $k$-tuples $(x_1\l x_k)$ of nonnegative integers such that $x_1 + \cdots + x_k\le n$.

\begin{thm}\label{vectorpath}
There is an algorithm as follows:
\begin{itemize}
\item{{\bf Input:}} A digraph $H$, and distinct vertices $s_0,t_0\in V(H)$; an integer $n\ge 0$; and for each edge $e$ of $H$, 
a member $l(e)$ of $K_n$.
\item{{\bf Output:}} The set of all minimal (under component-wise domination) vectors $l(P)$, over all paths $P$ of $H$ from $s_0$ to $t_0$; where for a path $P$ 
with edge set $\{e_1\l e_p\}$,
$l(P) = l(e_1)+\cdots + l(e_p)$.
\item{{\bf Running time:}} $O(n^{k}|V(H)||E(H)|)$.
\end{itemize}
\end{thm}
\Proof 
Let $Q_0(s_0)= \{(0\l 0)\}$, and let $Q_0(v) = \emptyset$ for every other vertex $v$ of $D$. Inductively, 
for $1\le i\le |V(H)|$, let $Q_i(v)$ be the set of minimal vectors in $K_n$ that either belong to $Q_{i-1}(v)$ or are
expressible in the form $l(e)+x$ for some edge $e = uv$ of $H$ and 
some $x\in Q_{i-1}(u)$.

Now here is an algorithm for the problem:
\begin{itemize}
\item For $i = 1\l |V(H)|$ in turn, compute $Q_i(v)$ for every $v\in V(H)$. 
\item Output $Q_{|V(H)|}(t_0)$.
\end{itemize}

It is easy to check that this output is correct, and we leave it to the reader. 
To compute $Q_i(v)$ at the $i$th step takes time $O(n^{k})d^-(v)$, where $d^-(v)$ is the in-degree of $v$ in $H$ (since $K_n$ has at most $(n+1)^k$ members),
and so the $i$th step in total takes time $O(n^{k}|E(H)|)$.
Thus the running time is $O(n^{k}|V(H)||E(H)|)$.~\bbox

Finally, we can give the main algorithm, \ref{mainchar}, which we restate.

\begin{thm}\label{mainchar2}
For all $d,k\ge 1$, there is an algorithm as follows:
\begin{itemize}
\item{{\bf Input:}} A problem instance $(G,s_1,t_1\l s_k,t_k)$ where $G$ is $d$-path-dominant.
\item{{\bf Output:}} The set of all key qualities of $(G,s_1,t_1 \l s_k,t_k)$.
\item{{\bf Running time:}} $O(n^t)$ where $t = 6k^2d(k+d) +13k$.
\end{itemize}
\end{thm}
\Proof Here is the algorithm. 
\begin{itemize}
\item Compute the $(k,m,c)$-tracker $H$, where 
\begin{eqnarray*}
c &=& ((k-1)d+k^2+2)(2k+1)k^2\\
m &=& ((k-1)d+k^2+2)d+1.
\end{eqnarray*}
\item For each edge $e = uv$ of $H$, define $l(e)$ as follows:
\begin{itemize}
\item if $u = s_0$ and $v = (L, X,Y)$ where $L = (M_1\l M_k)$, let $l(e) = (|V(M_1)|\l |V(M_k)|)$;
\item if $u = (L,X,Y)$ where $L = (M_1\l M_k)$, and $v = (L',X',Y')$ where $L' = (M_1'\l M_k')$, let $l(e) = (|V(M_1')\setminus V(M_1)|\l |V(M_k')\setminus V(M_k)|)$;
\item if $v = t_0$ let $l(e) = (0\l 0)$.
\end{itemize}
\item Run the algorithm of \ref{vectorpath} with input $H, s_0, t_0, l$.
\item Output its output.
\end{itemize}
To see its correctness, we must check that every key quality is in the output, and everything in the output is a key quality. 
We show first that every vector in the output is a quality. For let $x$ be in the output, and let $P$ be a path in $H$ from $s_0$ to $t_0$
with $l(P) = x$. By \ref{pathtolink}, $P$ traces some linkage $L = (P_1\l P_k)$ for $(G, s_1, t_1\l s_k, t_k)$; and so
$(|V(P_1)|, |V(P_2)|\l|V(P_k)|) = l(P)= x$. Hence $x$ is a quality. 

Next, we show that every key quality is in the output. For let $x$ be a key quality.
Let $L$ be an $x$-linkage for $(G, s_1, t_1\l s_k, t_k)$. By \ref{linktopath}, there is a path $P$ of $H$ from $s_0$ to $t_0$
tracing $L$; and hence $l(P) = x$ (where $l(P)$ is defined as in the statement of \ref{vectorpath}). 
Thus the output of \ref{vectorpath} contains a vector dominated by $x$. But $x$ does not dominate any other
quality, since it is key; and since every member of the output is  a quality, it follows that $x$ belongs to the output.

Third, we show that every member of the output is key. For let $x$ be in the output, and suppose it is not key. Hence $x$ dominates some other
quality, and hence dominates some other key quality $y$ say. Consequently $y$ is in the output. But no two members of the output dominate one another,
a contradiction. This proves that every member of the output is key, and so completes the proof that the output of the algorithm is as claimed.

Finally, for the running time: by \ref{countrails}, we can find all $(k,m,c)$-rails in time $O(n^{2km+1})$; and since there are
at most $O(n^{2km})$ of them (by \ref{countrails}), we can compute $H$ and the function $l$ in time $O(n^{4km})$. 
Then running \ref{vectorpath} takes time $O(n^{k}|V(H)|^3)$, and hence time at most $O(n^{6km+k})$.
Thus the total running time is $O(n^{6km+k})$.
Since $m = ((k-1)d+k^2+2)d+1$, the running time is $O(n^t)$ where 
$$t = 6k(k-1)d^2 + 6k(k^2+2)d + 7k = 6k^2d^2 + 6k^3d + 12 kd+ 7k - 6kd^2\le 6k^2d(k+d) +13k$$
as claimed. This proves \ref{mainchar2}.~\bbox


\begin{thebibliography}{99}

\bibitem{bangtho} J{\o}rgen Bang-Jensen and Carsten Thomassen , ``A polynomial algorithm for the 2-path problem for semicomplete
digraphs'', {\em SIAM Journal on Discrete Mathematics} 5 (1992), 366--376.
\bibitem{FHW} S.Fortune, J.Hopcroft and J.Wyllie, ``The directed subgraphs homeomorphism problem'',
{\em Theoret. Comput. Sci.} 10 (1980), 111--121.
\bibitem{FS} Alexandra Fradkin and Paul Seymour, ``Edge-disjoint paths in digraphs with bounded independence number'', 
{\em J. Combinatorial Theory, Ser.B}, in press, 2014 (DOI: 10.1016/j.jctb.2014.07.002).
\bibitem{RS} N. Robertson and P.D. Seymour, ``Graph minors. XIII. The disjoint paths problem'',
{\em J. Combinatorial Theory, Ser. B}, 63 (1995), 65--110.

\end{thebibliography}
\end{document}